\documentclass[a4paper,12pt]{article}
\usepackage[english]{babel}
\usepackage{amssymb,amsmath,amsthm,amscd}
\usepackage[a4paper]{geometry}
\usepackage[all,cmtip]{xy}
\usepackage{url}
\newcommand{\Z}{\mathbb{Z}}
\newcommand{\Q}{\mathbb{Q}}

\makeindex

\title{About the Proof of Lang's Height Conjecture}
\author{Benjamin WAGENER}

\begin{document}

\maketitle

\abstract{\textit{This informal document provides a concise overview of results related to the the recent paper} \cite{Wagener}.}

\medskip

Lang's height conjecture by itself was an important conjecture about the arithmetic of
elliptic curves. As it says, it involves the notion of height which is a kind of arithmetic-geometric
complexity and it makes the comparison between the height of rational points of elliptic curves
and the height of the curve itself. A comparison that is uniform with respect to the number field
on which the curve is defined which means that the comparison involves the height of rational
points on the curve, the height of the curve and a comparison term that depends only on the field
of definition of the curve.

As such is says roughly that "\textit{the complexity of a rational point of an elliptic curve is uniformly bounded from below by the complexity of the curve itself}".

\textbf{Important Note:} \textsl{During those about five years I worked intensively on this question, I was notably obsessed by the fact that there may exist an "idealistic" way to deal with it. I mean that if such result exists it should possibly be explained by a kind of "universal" argument that would probably take only a few pages. This is perhaps not false but not for the moment... We don't know enough to do this for now and probably still for very long. In contrast to this "Proof from the Book", the proof I propose may look very complicated to many, even possibly to the experts. Those that already practiced some Diophantine Mathematics know how much they are usually hard and tricky, those that practiced some Diophantine Geometry know this can be even more. From this point of view, the proof won't much chock the diophantine experts but it can chock even the Arithmetic Geometry expert. However you should know that the strategy of proof is quite natural. Remark that the paper itself is also especially complicated because, apart from the proof itself, I also did a hard job of optimization in order that the final result be as precise as it reasonably can be.}

\textsl{Note that I treated the question as "Problem-Solving" this is why the abstract of \cite{Wagener} says "focuses on the proof...". This paper has to be momentarily considered as such even if it can look very hard.}

\textbf{Remark about the strategy of proof:}
\begin{enumerate}
	\item \textsl{The strategy uses crucially the decomposition of the canonical height into local heights. This is natural as with respect to the meaning of the terms present in the discriminant and it is also natural with respect to the relation between each local height and the corresponding term in the discriminant.}
	\item \textsl{We uses a "Transcendental strategy", this is the kind of methods that was developed in the early 20th century when dealing with transcendental numbers. The reader may nicely note that this dates back from the 19th century with the proof of the transcendence of "$e$" by Hermite. From this point of view the method we use in \cite{Wagener} has been greatly Truth-tested and already provided various results.}
	\item \textsl{I didn't began this work completely from scratch. Marc Hindry told me that the conjecture was hard especially in a case when such a transcendental method was available. But there is something essential that explains why it was worthy to spend so much time so hardly trying to prove this: if one tries to apply an elementary strategy of "auxiliary polynomials" with the Baker method, we get already close to the result. This elementary, rough indeed, method fails because the estimates we can do with this are not precise enough but they are already quite.}
	\item \textsl{This is why it looked even at first that a complete proof of height conjecture was perhaps possible but we needed something very precise. This was fine, it was necessary to use Arakelov Theory and in spite of Baker Method, its (sharp) arakelovian equivalent: the slope method. It was nevertheless remaining a lot to do...}
	\item \textsl{In the paper we therefore use some of the most advanced knowledge and tools, as well from a Diophantine Mathematics point of view as from an Arithmetic Geometry point of view.}
\end{enumerate}

This conjecture appeared first in one of the many books of Serge Lang, Elliptic Curves/Diophantine Analysis,
dated from 1978 in the following form:

\textit{"... We select a point $P_1$ of minimum height $\neq 0$ (the height is the N\'eron-Tate height). We then let $P_2$ be a
point linearly independent of $P_1$ of height minimum $\geq 0$. We continue in this manner to obtain a maximum set of
linearly independent elements $P_1,\ldots, P_r$ called the} \textbf{successive minima}.

\textit{We suppose that A is defined by a Weierstrass equation}
\[y^2=x^3+ax+b,\]
\textit{and we assume for simplicity that $K=\mathbb{Q}$, so that we take $a,b\in\mathbb{Z}$....}

\textit{It seems a reasonable guess that uniformly for all such models of elliptic curves over $\mathbb{Z}$, one has}
\[\hat h(P_1)\gg \log|\Delta|.\]
\textit{It is difficult to guess how the next successive minima may behave. Do the ratios of the successive heights remain
bounded by a similar bound involving $|\Delta|$, or possibly of power of $\log|\Delta|$, or can they become much larger?..."}

Later in the same book, at the page 140, we can read an other conjecture:

\textit{"... In this connection, one could give a general setting for some of Demjanenko's ideas
[De 5], [De 6], and I would conjecture that on quasi-minimal model of an elliptic curve over a number
field $K$, the number of integral points is bounded by a number depending only on the rank of the
Mordell-Weil group (and of course K). However, to pursue these ideas would require knowing a lower
bound for the N\'eron-Tate height of point of infinite order, and such bounds are not known. Furthermore,
even in the special case considered in [De 6], I cannot understand his lemma 4 which seems to overlook
completely what are the essential difficulties in dealing with this type of problem."}

Indeed, in the work we did, we prove both conjectures, the second following height conjecture as said
by Serge Lang and we also prove an other theorem on abelian varieties that it also implies after some work.

We won't look at any detail here but we can see an elliptic curve as a curve $E$ defined by a Weierstrass equation
as above and the discriminant $\Delta$ that appears is a fundamental invariant of the curve that is directly
linked to the discriminant of the corresponding cubic equation. It is due to Joseph Silverman to have
seen that this conjecture can be restated with the Faltings height of the curve, $h_F(E)$,  replacing the
logarithm of the discriminant.

The Faltings height is generally defined for any abelian variety over a number field by means of Arakelov
theory. This is a fundamental scalar invariant defined as the Arakelov degree of the bundle of invariant differential forms of the N\'eron Model. It is therefore a fundamental numerical invariant associated to an abelian variety.

Such lower bound for the height of rational points is very important to know in Arithmetic Geometry and
Diophantine Geometry because it implies some important consequences. This lower bound is the "lowest complexity", something that has necessarily some intrinsic meaning, the fact that this "lowest(intrinsic) complexity" of rational points is directly equivalent to the "complexity" of the curve seems reasonable. What was surprising
in this conjecture is the uniformity of the inequality with respect to the number field of definition
but also, just knowing a lower bound would already have been interesting.

In our proof, we improved moreover the conjecture in two ways, first by proving that the result holds uniformly with respect
to the degree of the number field of definition and then by proving that moreover each term in the conjecture
is effectively computable. 

Moreover, the proof itself provides some bound on the torsion of elliptic curves, known
as Mazur-Merel theorem, which is also effectively computable.

The result we proved reads as follows:

\textbf{\textit{Theorem}} (Wagener (2017))

\textit{Let $d\geq 1$ be an integer,
there exist some constants $B_d$, $C_d$, $C_d'$, effectively computable, depending only on $d$ so that
for any elliptic curve
$E$ defined over a number field $K$ of degree $d$
and for any rational point $P\in E(K)$:}

\begin{enumerate}
\item \textit{either $P$ is a torsion point and there exists $n\leq B_d$ such that}
 \[[n]P=O_E.\]
\item \textit{or $P$ is of infinite order, and then, if $\Delta_{E/K}$ denotes the minimal discriminant ideal
of $E/K$ and $h_F(E/K)$ its Faltings height} \footnote{We don't suppose that it is the stable
Faltings  height.}, \textit{the following two inequalities are simultaneously true:}

\[\hat{h}(P)\geq C_d\log|N_{K/\mathbb{Q}}\Delta_{E/K}|,\]
and
\[\hat{h}(P)\geq C_d'\; h_F(E/K).\]
\end{enumerate}

Moreover the previous inequalities are satisfied with:
	\begin{gather*}
B_d=10^{16}(d\log(2d))^{1.54}\\
C_d=\frac{1}{10^{16}d^2(6094080+4.682*10^{15}d+883200d\log(48d))^2}\\
C_d'=\frac{1}{10^{12}d^2(6094080+4.682*10^{15}d+883200d\log(48d))^{2.08}}
\end{gather*}

\medskip

\textbf{Remarks:} \begin{itemize}
	\item \textsl{
The idea that proving the Height Conjecture could induce a bound on the torsion seems to originate from John Tate himself.}
	\item
\textsl{I didn't treated in the first versions of this paper the most elementary case in which $log|\Delta|=0$. Actually it was so trivial in comparison to the important cases I had to treat that I completely forgot... This is the case of everywhere good reduction. The former version the theorem is still true in this case but not interesting as we already know that the local height is positive. But actually we can quite easily say more in this case: we can easily prove that in this case the canonical height is bounded from below by a positive constant depending only on the degree. (For those that have had some look to the proof, in this case we only have to consider the possible issue of archimedean places... that anyway I treat in full generality in \cite{Wagener}...)}
	\item \textsl{Just a short comment about Serge Lang "seemingly reasonable guess" (see Serge Lang above). The proof of this "reasonable guess" is finally logically quite extreme. Apart from the "Transcendental method" that is logically already very tough, the geometric tools we use are very tough also but there is especially a vast necessary work of balancing and optimization. The proof conclude when I prove that the worst possible case is nevertheless also true. For this I have to evaluate a map on more than 16 Million points... It is funny to see what a "reasonable guess" of a Mathematician of the class of Serge Lang can be...}
	
\end{itemize}

In the way we proceed, Mazur-Merel's theorem appears complementary to our main result and
we find back this theorem by a quite easy argument when we turn to consider in the proof that the points
are torsion points.

In summary we obtain the following:

\begin{enumerate}
\item
A proof of Lang's height Conjecture.
\item
In the proof the uniformity appears as a uniformity with respect to the degree.
\item
We find back Mazur-Merel theorem in a completely new way and in complement to height theorem but, as has been long awaited, we provide a (small) polynomial bound.
\item
The result is effective in the sense that every term can be computed.
\item Our estimates are rather sharp.
\item
The bound for $B_d$ is of the form $B_d=10^{16}(d\log(2d))^{1.54}$ and it can be shown that a bound of this form is almost optimal. This bound is moreover polynomial (and smally polynomial), something that was conjectured since the first results (see below) on the torsion.
\item We can moreover note that thanks to the global geometric means we use, the Szpiro ratio (see below) appears nowhere. As such our work is completely independent of the Szpiro's Conjecture.
\end{enumerate}

The proof itself is quite intricate but make a very interesting use of Arithmetic Geometry
and Arakelov Geometry. There is absolutely no modular form, we make
an original arithmetic-geometric construction based on some diophantine argument in the
setting of the Slope Method that was developed in the 90s by Jean-Benoît Bost. This
geometric proof is very interesting by itself and we hope to open our way into more general
Arithmetic Geometric and Arakelov constructions that has been hoped to be possibly
performing since the beginning of Arakelov theory.

Apart from some tricks (this means "a lot" for those than never saw some diophantine mathematics/geometry), we can suggest why those results are rather sharp by the fact we work with global geometric concepts that tends to encapsulate in themselves the core of elliptic curves structure.
From this point of view we greatly benefit from Grothendieck point of view in arakelovian settings and from all the works of those that followed this view because of their great meaning.

The first result known on the torsion of elliptic curve is the following famous theorem. 

\textbf{Theorem} (Mazur (1977), \cite{Mazur1,Mazur2})

\textit{Let $E/\Q$ be an elliptic curve defined over $\Q$, then the possible groups of rational torsion points $E(\Q)_{\textrm{tors}}$ are given by,
	}
\[\left\{\begin{array}{l}\Z/m\Z\;\textrm{avec}\;m\in\{1,2,3,4,5,6,7,8,9,10,11,12\}\\\textrm{ou}\\
\Z/2m\Z\times\Z/2\Z\;\textrm{avec}\; m\in\{1,2,3,4\}\end{array}\right.\]

\medskip

Mazur proof made crucial uses of modular forms and the corresponding technics in details. However those technics are known for the moment only truly on the rational numbers $\mathbb{Q}$. From this point of view, the almost "exact" precision of this theorem that provides not only a bound but also the exact forms of the possible torsion groups, will remain for long an unreachable dream in the case arbitrary number fields. This would at least request a precise knowledge of modular/automorphic forms over arbitrary number fields that even in the case of elliptic curves we don't have...

This question related to the "complexity" of the ground field is more easy in the geometric context we use, perhaps simply because Grothendieck already made all the job for those mathematics to be not only valuable on any ground field but also remain consistently the same thanks to the global theoretical tools he(and others) developed.

\medskip

A primary result for the torsion for number field of arbitrary degree was given by:

\textbf{Theorem} (Flexor-Oesterl\'e (1990), \cite{flex})

\textit{Let $E/K$ be an elliptic curve defined over a number field $K$ of degree $d$.}
\begin{enumerate}
	\item \textit{If $E$ has bad reduction of additive type in at least one  place:}
	\[\textrm{Card}(E(K)_{\textrm{tors}})\leq 48 d.\]
	\item \textit{If $E$ has bad reduction of additive type in at least two places of different residual characteristics}
	\[\textrm{Card}(E(K)_{\textrm{tors}})\;\textrm{divides}\; 12,\]
	\textit{and if moreover one of this characteristic is greater or equal to $5$}
	\[\textrm{Card}(E(K)_{\textrm{tors}})\leq 4.\]
	\item \textit{If $E$ has good reduction in a place of residual characteristic $2$,}
	\[\textrm{Card}(E(K)_{\textrm{tors}})\leq 5\cdot 2^d.\]
\end{enumerate}

Next, Kamienny generalized the work of Mazur to quadratic number fields:
 
\textbf{Theorem} (Kamienny (1992), \cite{Kam})

\textit{Let $E/K$ be an elliptic curve defined over a number field $K$ of degree 2 over $\Q$. Then the group $E(K)_{\textrm{tors}}$
	is of finite order independent of $K$ and moreover this order is not divisible by any prime number $>13$.}

\medskip

The next big step was made by Merel that reworked the results of Mazur and Kamienny to obtain:

\medskip

\textbf{Theorem} (Merel (1996), \cite{Merel})

\textit{Let $d$ un entier $\geq 1$. There exists a positive number $B(d)$ such that for any elliptic curve  $E$, defined on a number field $K$
	of degree $d$ over $\Q$, any torsion point of $E(K)$ be of order $\leq B(d)$. In the case were the point is of $p-$torsion with $p$ prime
	one has:}
\[p\leq d^{3d^2}.\]

\medskip

This result is nevertheless not effective, it doesn't provide an effective formula for $B(d)$ but says that given the bound on primary torsion, results of Faltings and Frey say after Merel's work that this bound exists.

Note that Oesterl\'e, in an unpublished paper, improved the bound $d^{3d^2}$ to $(1+3^{\frac{d}{2}})^d$.

Note that Parent in \cite{Parent} obtain in 1999 a generalization of Oesterl\'e-Merel-type bounds to bounds on $p^n$-torsion points.

A much better bound was suggested by Hindry and Silverman but only in the simplest case of everywhere good reduction:

\medskip

\textbf{Theorem:} (Hindry-Silverman (1999), \cite{Hindry3})

\textit{Let $E/K$ be an elliptic curve on a number field $K$ with everywhere good reduction, then}

\[\textrm{Card}(E(K)_{\textrm{tors}})\leq 1977408 d\log\; d.\]

\medskip

It seems that Hindry and Silverman initiated some work on elliptic curves with complex multiplication to look for a global bound on torsion. They found in \cite{Hindry3} that this bound is  $>C^{te}\sqrt{d\log(\log(d))}$. Later Clark and Pollack \cite{Clarke1,Clarke2}, that for CM elliptic curves the optimal bound (for $d>>1$) on the cardinality of rational torsion is provided by:

\[\frac{e^{\gamma}\pi}{\sqrt{3}}d\log(\log(d)),\]

Other results were obtained by Breuer, Rattazzi, Masser and others.

\medskip

Thus our result implies that the optimal polynomial exponent, is between 1 (thanks to the last one) and 3.08 (thanks to our bound).

\medskip

\textbf{Secondary Remark:} In the paper \cite{Wagener} there is a secondary typo saying the optimal exponent is in between 1 and 1.54. Actually a bound on the order of torsion implies, this is well known, a bound on the cardinality of the rational torsion which is its square. Thus the proper upper bound on the polynomial exponent is $3.08$... not $1.54$...

\medskip

About Lang's conjecture, the first result was obtained by Silverman:

\medskip

\textbf{Theorem} (Silverman (1981), \cite{Joe6})

\textit{Fix a number fields $K$. There exists constants $c_1,c_2>0$ depending only on $K$ such that for any elliptic curve $E/K$ with integral $j-$invariant and for any point of infinite order}
 $P\in E(K)$, 
\[\hat{h}(P)>c_1\log|N_{K/\Q}\Delta_{E/K}|+c_2.\]

\medskip

A second result in the direction of Lang's conjecture was obtained by Hindry and Silverman, \cite{Hindry2} and comes from the idea, that seems to be due to Lang himself that the Szpiro Conjecture could be used to prove Lang Conjecture. Szpiro Conjecture, see for example \cite{Szpiro} Chapter 2,
suggests that for any number field $K$ and for all $\epsilon >0$, there exist only a finite number of elliptic curves
$E$ defined on $K$ for which
\[\sigma_{E/K}\geq 6+\epsilon,\]
where $\sigma_{E/K}$ is called the Szpiro ratio and is defined by
\[\sigma_{E/K}=\frac{\log |N_{K/\Q}\Delta_{E/K}|}{\log|N_{K/\Q}F_{E/K}|},\]
where $\Delta_{E/K}$ is the minimal discriminant ideal of $E/K$ and $F_{E/K}$ is its conductor ideal.

\medskip

\textbf{Theorem} (Hindry-Silverman (1988), \cite{Hindry2})

\textit{For any number field $K$ and for any elliptic curve $E$ defined on $K$,}

\[\hat{h}(P)\geq (20\sigma_{E/K})^{-8[K:\Q]}10^{-4\sigma_{E/K}} h_K(E),\]

This last "$h_K(E)$" being comparable to the Faltings height.

\medskip

Note that in the same article Hindry and Silverman established the truth of Lang's conjecture for fonction fields.

By an other method David \cite{Dav2} obtained in 1997 a lower bound of the following form:

\[\hat{h}(P)\geq c\frac{h}{[K:\Q]^3 \sigma_{E/K}^5}\left(1+\frac{\log\;[K:\Q]}{h}\right)^{-2},\]
where $h=\textrm{max}\{1,h(j_{E/K})\}$ o\`u $h(j_{E/K})$ is the absolute logarithmic Weil height of the modular invariant of $E/K$.

After those works it has been believed that a proof of Lang's Conjecture would require a proof of Szpiro's Conjecture. It is interesting to note that our demonstration is intrinsic enough in order that the Szpiro ratio doesn't appear.

\medskip
Knowing that Lang conjecture is true, we can easily deduce the following theorem by working back the proof of theorem 9.1 of \cite{Hindry2}. This is quite easy mainly after the Height's theorem and some previous work... But actually this is a very deep result also, Siegel's Theorem is not new but this form is effective and uniform. This is the second result that was conjectured by Serge Lang in his book (see above):

\medskip

\textbf{Theorem} (A Uniform  Siegel's Theorem)

\textit{For any number field $K$, there exists a constant $\mathcal{C}(K)$, effectively computable,
that depends only on $K$ such that for any finite set of places $S$ of $K$ that contains at least
all archimedean places and for any elliptic curve $E$ defined over $K$ given by a quasi-minimal
equation, for $\mathcal{O}_{K,S}-$ the ring
of $S-$integers of K, $|S|$ the cardinality of $S$ and $\textrm{rk}(E(K))$ the Mordell-Weil rank
of $E(K)$}:

\[\textrm{Card}\left(E(\mathcal{O}_{K,S})\right)\leq\mathcal{C}(K)^{|S|+1+\textrm{rk}(E(K))}.\]

\medskip

We can also deduce the following uniformity theorem by reconsidering the work of De Diego in \cite{Teresa} that was missing the Height's Theorem.

\medskip

\textbf{\textit{Theorem:}} (A Partial uniform  Faltings Theorem)

\textit{Let $d,d'\geq 1$ be some integers, $K$ a number field of degree $d$ and $K'$ some extension of $K$ of relative degree $d'=[K':K]$.}

 \textit{For all family of curves parametrized by some curve $T$, $\pi:\mathcal{C}\rightarrow T$,  if we suppose that $\mathcal{C}$ and
 $T$ are projectives on the number field $K$ and that the generic fiber $\mathcal{C}_\eta=\mathcal{C}\times_T\textrm{Spec}\;K(T)$
 of this family is smooth of genus $g\geq 2$. Denoting $T^0$ the open set of $T$ above which $\pi$ is smooth and defining
 moreover $\mathcal{J}\rightarrow T^0$ the family of jacobians associated to $\mathcal{C}$.}

\textit{Then, if there exist $g$ family of non-isotrivial elliptic curves $\mathcal{E}_1,\ldots,\mathcal{E}_g$ and a $T^0-$morphism $\alpha:\mathcal{J}\rightarrow \mathcal{E}_1\times\mathcal{E}_2\times\cdots\times \mathcal{E}_g$
such that for all $t\in T_0$, the fiber $\alpha_t:\mathcal{J}_t\rightarrow \mathcal{E}_{1t}\times\mathcal{E}_{2t}\times\cdots\times \mathcal{E}_{gt}$
is an isogeny.}

 \textit{There exists some constant $\mu=\mu_{d,d',\pi,\textrm{deg}(\alpha)}$, effectively computable, that depend only on $d$, on $d'$, on the family $\pi$ and on the degree of the isogeny $\alpha$ so that: }
	
\textit{For all	$t\in T^0(K')$}:

\[\textrm{Card}\left(\mathcal{C}_t(K')\right)\leq\mu^{(1+\textrm{rk}\mathcal{J}_t(K'))}.\]

This is quite a deep theorem also, especially because there is quite a good proportion of Jacobians that are isogenous to a product of elliptic curves... But for more we will need a "Height's Theorem" on Abelian varieties (or at least on Jacobians)... 
\bibliographystyle{alpha}
\bibliography{bibliothese-1}
\end{document}